\documentclass[11pt]{amsart}
\usepackage{geometry}                
\usepackage{graphicx}
\usepackage{amssymb}
\usepackage{epstopdf}
\title{On Optimal Offensive Strategies in Basketball}
\author{Ikjyot Singh Kohli}
\date{June 20, 2015}                                           

\begin{document}
\maketitle

\begin{abstract}
The purpose of this paper is to determine whether basketball teams who choose to employ an offensive strategy that involves predominantly shooting three point shots is stable and optimal. 
We employ a game-theoretical approach using techniques from dynamical systems theory to show that taking more three point shots to a point where an offensive strategy is dependent on predominantly shooting threes is not necessarily optimal, and depends on a combination of payoff constraints, where one can establish conditions via the global stability of equilibrium points in addition to Nash equilibria where a predominant two-point offensive strategy would be optimal as well. We perform a detailed fixed-points analysis to establish the local stability of a given offensive strategy. We finally prove the existence of Nash equilibria via global stability techniques via the monotonicity principle.
We believe that this work demonstrates that the concept that teams should attempt more three-point shots because a three-point shot is worth more than a two-point shot is therefore, a highly ambiguous statement.
\end{abstract}

\tableofcontents

\section{Introduction}
We are currently living in the age of analytics in professional sports, with a strong trend of their use developing in professional basketball. Indeed, perhaps, one of the most discussed results to come out of the analytics era thus far is the claim that teams should shoot as many three-point shots as possible, largely because, three-point shots are worth more than two-point shots, and this somehow is indicative of a very efficient offense. These ideas were mentioned for example by Alex Rucker \cite{rucker} who said ``When you ask coaches what's better between a 28 percent three-point shot and a 42 percent midrange shot, they'll say the 42 percent shot. And that's objectively false. It's wrong. \emph{If LeBron James just jacked a three on every single possession, that'd be an exceptionally good offense}. That's a conversation we've had with our coaching staff, and let's just say they don't support that approach.''  It was also claimed in the same article that ``The analytics team is unanimous, and rather emphatic, that every team should shoot more 3s including the Raptors and even the Rockets, who are on pace to break the NBA record for most 3-point attempts in a season.'' These assertions were repeated in \cite{rucker2}. In an article by John Schuhmann \cite{schuhmann}, it was claimed that ``It's simple math. A made three is worth 1.5 times a made two. So you don't have to be a great 3-point shooter to make those shots worth a lot more than a jumper from inside the arc. In fact, if you're not shooting a layup, you might as well be beyond the 3-point line. Last season, the league made 39.4 percent of shots between the restricted area and the arc, for a value of 0.79 points per shot. It made 36.0 percent of threes, for a value of 1.08 points per shot.''

The purpose of this paper is to determine whether basketball teams who choose to employ an offensive strategy that involves predominantly shooting three point shots is stable and optimal. Although this problem to the best of the author's knowledge has not been studied before in the literature, several studies that provide an in-depth quantitative analysis of various aspects of basketball games using statistical and game theoretical methods have been established in  \cite{Perse2009612}, \cite{Zhang20072190}, \cite{Kvam2006788}, \cite{Larkey1997596}, \cite{Popescu200681}, \cite{DeMelo2008695}, \cite{Clair20071163}, \cite{Dodge20121}, \cite{Zhou20081141} and references therein. 

We will employ a game-theoretical approach using techniques from dynamical systems theory to show that taking more three point shots to a point where an offensive strategy is dependent on predominantly shooting threes is not necessarily optimal, and depends on a combination of payoff constraints, where one can establish conditions via the global stability of equilibrium points in addition to Nash equilibria where a predominant two-point offensive strategy would be optimal as well. 

\section{The Dynamical Equations}
For our model, we consider two types of NBA teams. The first type are teams that employ two point shots as the predominant part of their offensive strategy, while the other type consists of teams that employ three-point shots as the predominant part of their offensive strategy. There are therefore two predominant strategies, which we will denote as $s_{1}, s_{2}$, such that we define
\begin{equation}
\mathbf{S} = \left\{s_{1}, s_{2}\right\}.
\end{equation}

We then let $n_{i}$ represent the number of teams using $s_{i}$, such that the total number of teams in the league is given by
\begin{equation}
N = \sum_{i =1}^{k} n_{i},
\end{equation}
which implies that the proportion of teams using strategy $s_{i}$ is given by
\begin{equation}
x_i = \frac{n_{i}}{N}.
\end{equation}

The state of the population of teams is then represented by $\mathbf{x} = (x_{1}, \ldots, x_{k})$. It can be shown \cite{webb} that the proportions of individuals using a certain strategy change in time according to the following dynamical system
\begin{equation}
\label{eq:dyn1}
\dot{x}_{i} = x_{i}\left[\pi(s_{i}, \mathbf{x}) - \bar{\pi}(\mathbf{x})\right],
\end{equation}
subject to
\begin{equation}
\label{eq:constr1}
\sum_{i =1}^{k} x_{i} = 1,
\end{equation}
where we have defined the average payoff function as
\begin{equation}
\label{eq:avpayoff}
\bar{\pi}(\mathbf{x}) = \sum_{i=1}^{k} x_{i} \pi(s_{i}, \mathbf{x}).
\end{equation}

Now, let $x_{1}$ represent the proportion of teams that predominantly shoot two-point shots, and let $x_{2}$ represent the proportion of teams that predominantly shoot three-point shots. Further, denoting the game action set to be $A = \left\{T, Th\right\}$, where $T$ represents a predominant two-point shot strategy, and $Th$ represents a predominant three-point shot strategy. As such, we assign the following payoffs:
\begin{equation}
\pi(T,T) = \alpha, \quad \pi(T,Th) = \beta, \quad \pi(Th, T) = \gamma, \quad \pi(Th,Th) = \delta.
\end{equation}

We therefore have that
\begin{equation}
\label{payoff1}
\pi(T,\mathbf{x}) = \alpha x_{1} + \beta x_{2}, \quad \pi(Th, \mathbf{x}) = \gamma x_{1} + \delta x_{2}.
\end{equation}
From \eqref{eq:avpayoff}, we further have that
\begin{equation}
\label{payoff2}
\bar{\pi}(\mathbf{x}) = x_{1} \left( \alpha x_{1} + \beta x_{2}\right) + x_{2} \left(\gamma x_{1} + \delta x_{2}\right).
\end{equation}

From Eq. \eqref{eq:dyn1} the dynamical system is then given by
\begin{eqnarray}
\label{eq:x1d}
\dot{x}_{1} &=& x_{1} \left\{ \left(\alpha x_{1} + \beta x_{2}   \right) - x_{1} \left( \alpha x_{1} + \beta x_{2}\right) - x_{2} \left(\gamma x_{1} + \delta x_{2}\right) \right\}, \\
\label{eq:x2d}
\dot{x}_{2} &=& x_{2} \left\{ \left( \gamma x_{1} + \delta x_{2}\right) -x_{1} \left( \alpha x_{1} + \beta x_{2}\right) - x_{2} \left(\gamma x_{1} + \delta x_{2}\right) \right\},
\end{eqnarray}
subject to the constraint
\begin{equation}
\label{eq:constr2}
x_{1} + x_{2} = 1.
\end{equation}
Indeed, because of the constraint \eqref{eq:constr2}, the dynamical system is actually one-dimensional, which we write in terms of $x_{1}$ as
\begin{equation}
\label{eq:x1d2}
\dot{x}_{1} = x_{1} \left(-1 + x_{1}\right) \left[\delta + \beta \left(-1 + x_{1}\right) - \delta x_{1} + \left(\gamma-\alpha\right)x_{1}\right].
\end{equation}

From Eq. \eqref{eq:x1d2}, we immediately notice some things of importance. First, we are able to deduce just from the form of the equation what the invariant sets are. Following \cite{ellis}, we note that for a dynamical system $\mathbf{x}' = \mathbf{f(x)} \in \mathbf{R^{n}}$ with flow $\phi_{t}$, if we define a $C^{1}$ function $Z: \mathbf{R}^{n} \to \mathbf{R}$ such that $Z' = \alpha Z$, where $\alpha: \mathbf{R}^{n} \to \mathbf{R}$, then, the subsets of $\mathbf{R}^{n}$ defined by $Z > 0, Z = 0$, and $Z < 0$ are invariant sets of the flow $\phi_{t}$. Applying this notion to Eq. \eqref{eq:x1d2}, one immediately sees that $x_1 > 0$, $x_1 = 0$, and $x_1 < 0$ are invariant sets of the corresponding flow.

Further, there also exists a symmetry such that $x_{1} \to -x_{1}$, which implies that without loss of generality, we can restrict our attention to $x_{1} \geq 0$.

\section{Fixed-Points Analysis}
With the dynamical system in hand, we are now in a position to perform a fixed-points analysis. There are precisely three fixed points, which are invariant manifolds and are given by:
\begin{equation}
P_{1}: x_{1}^{*} = 0, \quad P_{2}: x_{1}^{*} = 1, \quad P_{3}: x_{1}^{*} = \frac{\beta - \delta}{-\alpha + \beta - \delta + \gamma}.
\end{equation}
Note that, $P_{3}$ actually contains $P_{1}$ and $P_{2}$ as special cases. Namely, when $\beta = \delta$, $P_{3} = 0 = P_{1}$, and when $\alpha = \gamma$, $P_{3} = 1 = P_{2}$. We will therefore just analyze, the stability of $P_{3}$.
$P_{3} = 0$ represents a state of the population where all teams predominantly shoot three-point shots. Similarly, $P_{3} = 1$ represents a state of the population where all teams predominantly shoot two-point shots, We additionally restrict
\begin{equation}
0 \leq P_{3} \leq 1 \Rightarrow 0 \leq \frac{\beta - \delta}{-\alpha + \beta - \delta + \gamma} \leq 1, 
\end{equation}
which implies the following conditions on the payoffs:
\begin{equation}
\left[\delta < \beta \cap \gamma \leq \alpha   \right] \cup \left[\delta = \beta \cap \left(\gamma < \alpha \cup \gamma > \alpha  \right) \right] \cup \left[\delta > \beta \cap \gamma \leq \alpha \right].
\end{equation}

With respect to a stability analysis of $P_{3}$, we note the following. 

%
The point $P_{3}$ is a:
\begin{itemize}
	\item Local sink if: $\{\delta < \beta\} \cap \{\gamma > \alpha\}$,
	\item Source if: $\{\delta > \beta\} \cap \{\gamma < \alpha\}$,
	\item Saddle: if: $\{\delta = \beta \} \cap (\gamma < \alpha -\beta + \delta \cup \gamma > \alpha - \beta + \delta)$, or $(\{\delta < \beta\} \cup \{\delta > \beta\}) \cap \gamma = \frac{\alpha \delta - \alpha \beta}{\delta - \beta}$.
\end{itemize}

Further, the system exhibits some bifurcations as well. In the neigbourhood of $P_{3} = 0$, the linearized system takes the form
\begin{equation}
x_{1}' = \beta - \delta.
\end{equation}
Therefore, $P_{3} = 0$ destabilizes the system at $\beta = \delta$. Similarly, $P_{3} = 1$ destabilizes the system at $\gamma = \alpha$. Therefore, bifurcations of the system occur on the lines $\gamma = \alpha$ and $\beta = \delta$ in the four-dimensional parameter space.

\section{Global Stability and The Existence of Nash Equilibria}
With the preceding fixed-points analysis completed, we are now interested in determining global stability conditions. The main motivation is to determine the existence of any Nash equilibria that occur for this game via the following theorem \cite{webb}: If $\mathbf{x}^{*}$ is an asymptotically stable fixed point, then the symmetric strategy pair $[\sigma^{*}, \sigma^{*}]$, with $\sigma^{*} = \mathbf{x}^*$ is a Nash equilibrium.

We will primarily make use of the monotonicity principle, which says \cite{ellis} let $\phi_{t}$ be a flow on $\mathbb{R}^{n}$ with $S$ an invariant set. Let $Z: S \to \mathbb{R}$ be a $C^{1}$ function whose range is the interval $(a,b)$, where $a \in \mathbb{R} \cup \{-\infty\}, b \in \mathbb{R} \cup \{\infty\}$, and $a < b$. If $Z$ is decreasing on orbits in $S$, then for all $\mathbf{x} \in S$,
\begin{equation*}
\omega(\mathbf{x}) \subseteq \left\{\mathbf{s} \in \partial S | \lim_{\mathbf{y} \to \mathbf{s}} Z(\mathbf{y}) \neq \mathbf{b}\right\},
\end{equation*}
\begin{equation*}
\alpha(\mathbf{x}) \subseteq \left\{\mathbf{s} \in \partial S | \lim_{\mathbf{y} \to \mathbf{s}} Z(\mathbf{y}) \neq \mathbf{a}\right\}.
\end{equation*}

Consider the function
\begin{equation}
Z_{1} = \log \left(-1 + x_{1}\right).
\end{equation}
Then, we have that
\begin{equation}
\dot{Z}_{1}= x_{1} \left[\delta + \beta \left(-1 + x_{1}\right) - \delta x_{1} + x_{1} \left(\gamma - \alpha\right)\right].
\end{equation}

For the invariant set $S_1 = \{0 < x_{1} < 1\}$, we have that $\partial S_{1} = \{x_{1} = 0\} \cup \{x_{1} = 1\}$. One can then immediately see that in $S_{1}$,
\begin{equation}
\dot{Z}_{1} < 0 \Leftrightarrow \left\{\beta > \delta\right\} \cap \left\{\alpha \geq \gamma\right\}.
\end{equation}

Therefore, by the monotonicity principle, 
\begin{equation}
\omega(\mathbf{x}) \subseteq \left\{\mathbf{x}: x_{1} = 1 \right\}.
\end{equation}
Note that the conditions $\beta > \delta$ and $\alpha \geq \gamma$ correspond to $P_{3}$ above. In particular, for $\alpha = \gamma$, $P_{3} = 1$, which implies that $x_{1}^{*} = 1$ is globally stable. Therefore, under these conditions, the symmetric strategy $[1,1]$ is a Nash equilibrium. 

Now, consider the function
\begin{equation}
Z_{2} = \log \left(x_{1}\right).
\end{equation}
We can therefore see that
\begin{equation}
\dot{Z}_{2} = \left[-1 + x_{1}\right] \left[\delta + \beta\left(-1+x_{1}\right) - \delta x_{1} + \left(-\alpha + \gamma\right) x_{1}\right].
\end{equation}
Clearly, $\dot{Z}_{2} < 0$ in $S_{1}$ if for example $\beta = \delta$ and $\alpha < \gamma$. Then, by the monotonicity principle, we obtain that
\begin{equation}
\omega(\mathbf{x}) \subseteq \left\{\mathbf{x}: x_{1} = 0 \right\}.
\end{equation}
Note that the conditions $\beta = \delta$ and $\alpha < \gamma$ correspond to $P_{3}$ above. In particular, for $\beta = \delta$, $P_{3} = 0$, which implies that $x_{1}^{*} = 0$ is globally stable. Therefore, under these conditions, the symmetric strategy $[0,0]$ is a Nash equilibrium.

In summary, we have just shown that for the specific case where $\beta > \delta$ and $\alpha = \gamma$, the strategy $[1,1]$ is a Nash equilibrium. On the other hand, for the specific case where $\beta = \delta$ and $\alpha < \gamma$, the strategy $[0,0]$ is a Nash equilibrium. 

\section{Discussion}
In the previous section which describes global results, we first concluded that for the case where $\beta > \delta$ and $\alpha = \gamma$, the strategy $[1,1]$ is a Nash equilibrium. The relevance of this is as follows. The condition on the payoffs thus requires that
\begin{equation}
\label{eq:presult1}
\pi(T,T) = \pi(Th,T), \quad \pi(T,Th) > \pi(Th,Th).
\end{equation}
That is, given the strategy adopted by the other team, neither team could increase their payoff by adopting another strategy if and only if the condition in \eqref{eq:presult1} is satisfied. Given these conditions, if one team has a predominant two-point strategy, it would be the other team's best response to also use a predominant two-point strategy.

We also concluded that for the case where $\beta = \delta$ and $\alpha < \gamma$, the strategy $[0,0]$ is a Nash equilibrium. The relevance of this is as follows. The condition on the payoffs thus requires that
\begin{equation}
\label{eq:presult2}
\pi(T,Th) = \pi(Th,Th), \quad \pi(T,T) < \pi(Th,T).
\end{equation}
That is, given the strategy adopted by the other team, neither team could increase their payoff by adopting another strategy if and only if the condition in \eqref{eq:presult2} is satisfied. Given these conditions, if one team has a predominant three-point strategy, it would be the other team's best response to also use a predominant three-point strategy.

Further, we also showed that $x_{1} = 1$ is globally stable under the conditions in \eqref{eq:presult1}. That is, if these conditions hold, every team in the NBA will eventually adopt an offensive strategy predominantly consisting of two-point shots. The conditions in \eqref{eq:presult2} were shown to imply that the point $x_{1} = 0$ is globally stable. This means that if these conditions now hold, every team in the NBA will eventually adopt an offensive strategy predominantly consisting of three-point shots. 

We also provided through a careful stability analysis of the fixed points criteria for the local stability of strategies. For example, we showed that a predominant three-point strategy is locally stable if $\pi(T,Th) - \pi(Th,Th) < 0$, while it is unstable if $\pi(T,Th) - \pi(Th,Th) \geq 0$. In addition, a predominant two-point strategy  was found to be locally stable when $\pi(Th,T) - \pi(T,T) < 0$, and unstable when $\pi(Th,T) - \pi(T,T) \geq 0$. 

There is also they key point of which one of these strategies has the highest probability of being executed. From \cite{webb}, we know that
\begin{equation}
\pi(\sigma,\mathbf{x}) = \sum_{s \in \mathbf{S}} \sum_{s' \in \mathbf{S}} p(s) x(s') \pi(s,s').
\end{equation}
That is, the payoff to a team using strategy $\sigma$ in a league with profile $\mathbf{x}$ is proportional to the probability of this team using strategy $s \in \mathbf{S}$. We therefore see that a team's optimal strategy would be that for which they could maximize their payoff, that is, for which $p(s)$ is a maximum, while keeping in mind the strategy of the other team, hence, the existence of Nash equilibria. 

Hopefully, this work also shows that the concept that teams should attempt more three-point shots because a three-point shot is worth more than a two-point shot is a highly ambiguous statement. In actuality, one needs to analyze what offensive strategy is optimal which is constrained by a particular set of payoffs. 

\newpage
\bibliographystyle{ieeetr}
\bibliography{sources}

\end{document}